\documentclass[preprint,12pt]{elsarticle}

\usepackage{amssymb}
\usepackage{amsmath}
\usepackage{amsthm}

\journal{Topology and its Applications}
\usepackage[T1]{fontenc}
\usepackage{lmodern}
\usepackage{microtype}
\usepackage{needspace}
\usepackage{placeins}
\usepackage{tikz}
\usepackage[hidelinks]{hyperref}
\newtheorem{theorem}{Theorem}[section]
\newtheorem{proposition}[theorem]{Proposition}
\newtheorem{lemma}[theorem]{Lemma}
\newtheorem{corollary}[theorem]{Corollary}
\theoremstyle{definition}
\newtheorem{definition}[theorem]{Definition}
\newtheorem{remark}[theorem]{Remark}
\newtheorem{example}[theorem]{Example}
\newcommand{\Afin}{\mathcal A_{\mathrm{fin}}}
\newcommand{\cl}{\operatorname{cl}}
\newcommand{\RO}{\operatorname{RO}}
\newcommand{\id}{\operatorname{id}}
\newcommand{\Int}{\operatorname{int}}
\newcommand{\Is}{\operatorname{Is}}
\newcommand{\Idl}{\operatorname{Idl}}

\usetikzlibrary{positioning}
\tikzset{
  vertex/.style={
    circle,
    fill=black,
    inner sep=0pt,
    minimum size=4pt
  },
  edge/.style={line width=.7pt},
  limit/.style={edge,densely dashed},
  every label/.style={font=\small}
}

\begin{document}

\begin{frontmatter}

\title{Meet-Continuity of Scott Function Spaces\footnote{The work was supported by the National Natural Science Foundation of China(Grants No. 12571495, 12231007)}} %% Article title

\author{Huanhuan He}
\author{Qingguo Li}
\author{Chengyu Zhou}
\author{Xiangnan Zhou\footnote{Corresponding author. xnzhou@hnu.edu.cn}}
\affiliation{addressline={School of Mathematics},
organization={Hunan University}, 
city={Changsha},
postcode={410082}, 
state={Hunan},
country={P.R.China}}

\begin{abstract}
We investigate the meet-continuity of the dcpo $[X\to L]$ of Scott-continuous maps from a poset $X$ to a meet-continuous dcpo $L$, without assuming that $L$ has binary meets. It is proved that $[X\to L]$ is meet-continuous whenever $X$ is an RB-domain. A counterexample shows that this result does not extend to all FS-domains. On the other hand, for a fixed meet-continuous dcpo $L$, the RB-domain hypothesis need not be necessary. Finally, we establish a characterization of RB-domains in terms of the meet-continuity of Scott function spaces: a continuous domain $X$ is an RB-domain if and only if $[X\to L]$ is meet-continuous for every meet-continuous dcpo $L$.
\end{abstract}

\begin{keyword}
Meet-continuity \sep Scott function space \sep RB-domain \sep
FS-domain \sep ordinal sum
\end{keyword}

\end{frontmatter}

\section{Introduction}\label{sec:introduction}
Scott function spaces play an important role in domain theory and denotational semantics. In denotational semantics, higher-order types are usually interpreted by spaces of Scott-continuous maps. For dcpos $X$ and $L$, we write $[X\to L]$ for the dcpo of Scott-continuous maps from $X$ to $L$, ordered pointwise. A basic question is to determine under what conditions on $X$ the function space $[X\to L]$ inherits a property of $L$, and conversely, under what conditions on $L$ it inherits a property of $X$ \cite{AJ,Gierz,JungCCC}.  

Several results are known for continuity and quasicontinuity of function spaces. Xi, Xu, and Lawson \cite{XiXuLawson2016} proved that the function space from a coherent continuous domain to a pointed RB-domain is again an RB-domain. Liu et al. \cite{LiuEtAl2020} showed that $[X\to L]$ is quasicontinuous when $X$ is a bounded complete algebraic domain and $L$ is a quasicontinuous domain with property $M^*$. Miao et al. \cite{MiaoEtAl} proved that $[X\to L]$ is quasicontinuous when $X$ is an RB-domain and $L$ is a quasicontinuous domain with property $M^*$. They also showed that $[X\to L]$ is continuous whenever $X$ is an RB-domain and $L$ is a continuous domain. More recently, Chen \cite{ChenFunction} showed that a continuous domain $X$ is an RB-domain if and only if $[X\to L]$ is continuous for every continuous domain $L$.

Meet-continuity connects quasicontinuous domains with continuous domains: a dcpo is continuous if and only if it is both quasicontinuous and meet-continuous\cite{KouLiuLuo}. It also plays an important role in the study of function spaces and Cartesian closed categories of domains \cite{JiaThesis,JiaJungLiDichotomy}. In this paper, we focus on the meet-continuity of Scott function spaces. Meet-continuous dcpos were introduced and studied by Kou, Liu, and Luo \cite{KouLiuLuo}, and further structural and topological results were obtained by Jia \cite{JiaThesis}. If $L$ is a meet-continuous dcpo with binary meets, then $[X\to L]$ is meet-continuous for every dcpo $X$ \cite[Lemma~5.5]{KouLiuLuo}. Jia, Jung, and Li \cite{JiaJungLiDichotomy}  proved that a dcpo $L$ is meet-continuous if its self-function space $[L\to L]$ is core-compact. They also showed that, if $X$ is a locally compact sober dcpo and $L$ is a pointed bicomplete sober dcpo, then meet-continuity of $[X\to L]$ implies that either $X$ is coherent or $L$ is an $L$-dcpo.

The assumption that $L$ has binary meets allows common lower bounds of functions to be constructed pointwise. Without this assumption, common lower bounds in $L$ need not have a greatest element, and lower bounds chosen at different points of $X$ may not fit together to form a Scott-continuous function. It is therefore not immediate how the meet-continuity of $L$ can be used to obtain meet-continuity of $[X\to L]$. The aim of this paper is to look for conditions on $X$ under which this preservation holds without assuming binary meets in $L$.

Inspired by the work of Miao et al. \cite{MiaoEtAl} on quasicontinuous function spaces, we prove that $[X\to L]$ is meet-continuous whenever $X$ is an RB-domain and $L$ is a meet-continuous dcpo. We first establish the result for finite posets and then extend it to bifinite domains using finite-image approximations. The result for RB-domains follows by Scott-continuous retractions.
It is also showed that meet-continuity of a function space does not imply that its Isbell and Scott topologies coincide. Using a result of Ern\'e and Gatzke \cite{ErneGatzke} on the regular-open algebra $\RO(\mathbb R)$, we give an example with a two-element discrete poset.

However, the preceding preservation result does not extend to all FS-domains. The closed-disk domain is an FS-domain \cite{AJ,Lawson2008} that it is
not an RB-domain \cite{ChenKouLyu}. Further examples separating FS-domains from RB-domains were obtained by Chen \cite{ChenCone} and Abbadini \cite{Abbadini}. Shen et al. \cite{ShenEtAl} constructed an FS-domain and an algebraic domain whose function space is not continuous. Using an algebraic cover construction adapted from their work and Chen's obstruction to finite-valued monotone approximation \cite{ChenCone}, we construct $\omega$-algebraic domains $L$ for which $[X\to L]$ is not meet-continuous, where $X$ is either the closed-disk domain or a domain associated with a non-simplicial proper cone.

The RB-domain condition, however, need not be necessary for a particular meet-continuous dcpo. This raises the question of whether the RB-domain condition becomes necessary when preservation is required for every meet-continuous dcpo.

For continuous domains, we answer this question affirmatively. In view of Chen's characterization \cite{ChenFunction}, it is natural to ask whether continuity of the function spaces can be weakened to meet-continuity. We show that a continuous domain $X$ is an RB-domain if and only if $[X\to L]$ is meet-continuous for every meet-continuous dcpo $L$. In fact, it is enough to require $[X\to A]$ to be meet-continuous for every algebraic domain $A$. For an algebraic domain $X$, these conditions are also equivalent to bifiniteness; compare the function-space characterization of bifinite domains by Luan and Li \cite{LuanLi2025}.

The remainder of this paper is organized as follows. Section 2 recalls some basic notions and results. In Section 3, we prove that the function space from an RB-domain to a meet-continuous dcpo is meet-continuous. Section 4 gives a counterexample showing that this result does not extend to all FS-domains. In Section 5, we study function spaces from ordinal sums of dcpos and shows that, for a continuous domain $X$, if $[X\to L]$ is meet-continuous for \emph{every} meet-continuous dcpo $L$, then $X$ is an RB-domain.

\section{Preliminaries}\label{sec:preliminaries}
Now we quickly recall some basic definitions and properties in domain theory. Readers can see \cite{AJ,Gierz} for more details. 

Let $P$ be a poset. For $A\subseteq P$, write $\mathord{\downarrow}A=\{x\in P:x\leq a\text{ for some }a\in A\}$ and define $\mathord{\uparrow}A$ dually; write $\mathord{\downarrow}a$ and $\mathord{\uparrow}a$ for the singleton. A subset $D\subseteq P$ is \emph{directed} if every finite subset of $D$ has an upper bound in $D$. A \emph{dcpo} is a poset in which every directed subset has a supremum. A dcpo with a least element is called \emph{pointed}. 

For a dcpo $P$, an upper set $U$ is \emph{Scott open} if $\bigvee D\in U$ for a directed set $D$ implies $D\cap U\neq\varnothing$. The \emph{Scott topology} consists of all Scott open sets, which is denoted by $\sigma(P)$. We write $\cl_{\sigma(P)}A$ for the Scott closure of $A$ and $\mathcal O(X)$ for the lattice of Scott-open subsets of $X$. The \emph{lower topology} of $P$ has subbasic open sets of the form $P\setminus{\mathord\uparrow x}$. The \emph{Lawson topology} of $P$ is a common refinement of Scott topology and lower topology. A dcpo is \emph{Scott compact} or \emph{Lawson compact} if it is compact in the corresponding topology. It is \emph{coherent} if the intersection of any two Scott compact subsets is Scott compact.

For $x,y$ in a dcpo $P$, we write $x\ll y$ if, whenever a directed set $D$ satisfies $y\leq\bigvee D$, there exists some $d\in D$ such that $x\leq d$. Define $\mathord{\Uparrow}x=\{y\in P:x\ll y\}$ and
$\mathord{\Downarrow}x=\{y\in P:y\ll x\}$. The set $\mathord{\Uparrow}x$ is Scott open whenever $P$ is continuous. An element $x$ is \emph{compact} if $x\ll x$, and the set of compact elements of $P$ is denoted by $K(P)$. The dcpo $P$ is \emph{continuous} if $x$ is the directed supremum of $\mathord{\Downarrow}x$ for every $x\in P$ and it is \emph{algebraic} if $x$ is the directed supremum of $\mathord{\downarrow}x\cap K(P)$. A \emph{domain} is a continuous dcpo, and an algebraic domain is called an \emph{algebraic domain}. A basis of a domain $X$ is a subset $B$ such that $\mathord\Downarrow x\cap B$ is directed with the supremum $x$ for every $x\in X$. The dcpo is countably based if it has a countable basis, and a countably based algebraic domain is called \emph{$\omega$-algebraic}.

An \emph{ideal} of a poset $P$ is a nonempty directed lower subset. The ideals of $P$, ordered by inclusion, form an algebraic domain $\Idl(P)$. Its directed suprema are unions, and its compact elements are exactly the principal ideals.

A map between posets is \emph{monotone} if it preserves the order. A map between dcpos is \emph{Scott-continuous} if it is monotone and preserves directed suprema. For dcpos $X$ and $L$, $[X\to L]$ denotes the dcpo of all Scott-continuous maps from $X$ to $L$, ordered pointwise. Directed suprema in $[X\to L]$ are computed pointwise. Consequently, for each $x\in X$, the evaluation map $\operatorname{eV}_x:[X\to L]\to L$, defined by
$\operatorname{eV}_x(f)=f(x)$, is Scott-continuous. A dcpo $X$ is a \emph{Scott-continuous retract} of a dcpo $Y$ if there are Scott-continuous maps $s:X\to Y$ and $r:Y\to X$ with $r\circ s=\id_X$.

\begin{definition}
Let $X$ be a dcpo. A \emph{deflation} is a Scott-continuous map $q:X\to X$ with finite image and $q\le\id_X$. A map $f$ is idempotent if $f\circ f=f$. An \emph{approximating identity} is a directed subset $\mathcal{D}\subseteq [X\to X]$ with $\bigvee\mathcal{D}=\id_X$. A Scott-continuous map $f:X\to X$ is \emph{finitely separated} if there is a finite set $M\subseteq X$ such that, for each $x\in X$, some $m\in M$ satisfies $f(x)\leq m\leq x$.    
\end{definition}

\begin{definition}
An \emph{FS-domain} is a dcpo with an approximating identity of finitely separated maps. A domain is an \emph{RB-domain} if it has a directed family of deflations with pointwise supremum $\id_X$. A dcpo is a \emph{bifinite domain} if it has an approximating identity consisting of idempotent deflations.
\end{definition}

A dcpo is an RB-domain if and only if it is a Scott-continuous retract of a bifinite domain. Every RB-domain is an FS-domain, but the converse does not hold \cite{ChenKouLyu}.

\begin{definition}
A dcpo $P$ is \emph{meet-continuous} if, for every directed set $D\subseteq P$ and every $x\leq\bigvee D$, $x\in\cl_{\sigma(P)}\bigl(\mathord{\downarrow}x\cap\mathord{\downarrow}D\bigr)$.
\end{definition}

A meet-continuous semi-lattice is a meet-continuous dcpo which has binary meets, i.e., a dcpo semi-lattice $P$ satisfies $x\wedge\bigvee D=\bigvee_{d\in D}(x\wedge d)$ for every $x\in P$ and every directed subset $D\subseteq P$.
Every domain is meet-continuous
\cite[Proposition~2.2]{KouLiuLuo}. Meet-continuity is preserved by Scott-continuous retractions \cite[Proposition~3.1.7]{JiaThesis}.

\begin{proposition}\cite[Lemma~5.5]{KouLiuLuo}\label{prop:pointwise}
Let $X$ be a dcpo and let $L$ be a meet-continuous semi-lattice.
Then $[X\to L]$ is a meet-continuous semi-lattice.
\end{proposition}

\section{Function spaces from RB-domains to meet-continuous dcpos}\label{sec:positive}

In this section, we try to look for the conditions of $X$, under which the function space $[X\to L]$ inherits the meet-continuity of the dcpo $L$. It is well known that every bifinite domain is a bilimit of finite posets. Following the approach of Miao et al.~\cite{MiaoEtAl}, we first investigate the function space from a finite poset to a meet-continuous dcpo. 

\begin{definition}\label{def:order-compatible}
Let $P$ and $Q$ be dcpos. A relation $R \subseteq P\times Q$ is
\emph{order-compatible} if
\[
u'\leq u,\quad uR v,\quad v\leq v'
\quad\Longrightarrow\quad u'R v'.
\]
\end{definition}

Such a relation is downward closed in its first variable and
upward closed in its second. Ordering $R$ coordinatewise, if
$D\subseteq P$ is directed with $dR v$ for every $d\in D$,
then $\{(d,v):d\in D\}$ is directed in $R$. 

\begin{lemma}\label{lem}
Let $P$ and $Q$ be meet-continuous dcpos, and let $R\subseteq P\times Q$ be a nonempty order-compatible relation. Suppose that directed suprema in $R$ exist and are computed coordinatewise. Assume that, for every directed set $E\subseteq Q$ and every $p\in P$ satisfying $pR\bigvee E$,
$$
p\in\cl_{\sigma(P)}
\{a\in P: a\leq p\text{ and }aRe\text{ for some }e\in E\}.
$$
Then $R$ is a meet-continuous dcpo.
\end{lemma}

\begin{proof}
We regard $R$ as a dcpo with the coordinatewise order. Let
$D=\{(p_i,q_i):i\in I\}\subseteq R
$
be directed, where $I$ is a directed index set, and suppose that
$
(p_0,q_0)\leq\bigvee_{i\in I}(p_i,q_i).
$
Take
$
C=\cl_{\sigma(R)}
\bigl(\mathord\downarrow_R(p_0,q_0)\cap\mathord\downarrow_R D\bigr).
$
To prove the meet continuity, we only need to show that $(p_0,q_0)\in C$.

For each $q\in Q$, define
$
C_q=\{p\in P:(p,q)\in C\}.
$
We first show that $C_q$ is Scott closed in $P$. If $p'\leq p$ and $p\in C_q$, then $pRq$ and order compatibility give $p'Rq$. Hence $(p',q)\in R$, and  $(p',q)\in C$ since $C$ is a lower set. Thus $p'\in C_q$. Hence $C_q$ is a lower set.

Now let $F\subseteq C_q$ be directed. Then
$\{(p,q):p\in F\}$
is a directed subset of $C$. Its supremum in $R$ is $(\bigvee F,q)$, because directed suprema in $R$ are computed coordinatewise. Since $C$ is Scott closed, $(\bigvee F,q)\in C$, and therefore $\bigvee F\in C_q$.

Define
\[
T_{p_0}=\{q\in Q: \text{for every }p\leq p_0,\
pRq\text{ implies }(p,q)\in C\}.
\]
We claim that $T_{p_0}$ is Scott closed in $Q$.

To show that $T_{p_0}$ is a lower set, let $q'\leq q\in T_{p_0}$ and suppose that $p\leq p_0$ and $pRq'$. By order compatibility, $pRq$. Since $q\in T_{p_0}$, we have $(p,q)\in C$. Thus $(p,q')\in C$. Hence $q'\in T_{p_0}$.

Next, let $E\subseteq T_{p_0}$ be directed. Fix $p\leq p_0$ such that $pR\bigvee E$. By assumption,
\[
p\in\cl_{\sigma(P)}\{a\in P: a\leq p\text{ and }aRe\text{ for some }e\in E\},
\]
Let $b\in \{a\in P: a\leq p\text{ and }aRe\text{ for some }e\in E\}$, and choose $e_0\in E$ such that $bRe_0$. For every $e\in E\cap\uparrow e_0$, order compatibility gives $bRe$. Moreover, $b\leq p\leq p_0$ and $e\in T_{p_0}$, so $(b,e)\in C$.

The set $E\cap\uparrow e_0$ is directed and cofinal in $E$. Consequently,
$
\bigvee_{e\in E\cap\uparrow e_0}(b,e)=(b,\bigvee E)
$
in $R$. Since $C$ is Scott closed, $(b,\bigvee E)\in C$.  As $C_{\bigvee E}$ is Scott closed, we obtain $p\in C_{\bigvee E}$, that is, $(p,\bigvee E)\in C$. Therefore $\bigvee E\in T_{p_0}$, and $ T_{p_0}$ is Scott closed.

We next show that
$
\mathord\downarrow_Q q_0\cap\mathord\downarrow_Q\{q_i: i\in I\}\subseteq  T_{p_0}.
$
Let $q\leq q_0,q_i$ for some $i\in I$. To prove that $q\in  T_{p_0}$, take any $p\leq p_0$ satisfying $pRq$. Since
$
p\leq p_0\leq\bigvee_{j\in I}p_j,
$
the meet-continuity of $P$ gives
$$
p\in\cl_{\sigma(P)}
\bigl(\mathord\downarrow_P p\cap\mathord\downarrow_P\{p_j: j\in I\}\bigr).
$$
Suppose that $a\leq p,p_j$ for some $j\in I$. Choose $k\in I$ above both $i$ and $j$. Since $a\leq p$ and $pRq$, order compatibility gives $aRq$. Moreover,
$
a\leq p\leq p_0, a\leq p_j\leq p_k,
$
and
$
q\leq q_0, q\leq q_i\leq q_k.
$
Thus
$
(a,q)\leq(p_0,q_0),(p_k,q_k),
$
so $(a,q)\in C$. Hence
$
\mathord\downarrow_P p\cap\mathord\downarrow_P\{p_j: j\in I\}\subseteq C_q.
$
Since $C_q$ is Scott closed, $p\in C_q$. Therefore $(p,q)\in C$, $q\in T_{p_0}$.

Finally, $q_0\leq\bigvee_i q_i$, so meet-continuity of $Q$ yields
$
q_0\in\cl_{\sigma(Q)}
\bigl(\mathord\downarrow_Q q_0\cap
\mathord\downarrow_Q\{q_i: i\in I\}\bigr).
$
The set inside the closure is contained in $T_{p_0}$, and $T_{p_0}$ is Scott closed. Hence $q_0\in T_{p_0}$. Since $(p_0,q_0)\in R$, we have $p_0Rq_0$. By the definition of $T_{p_0}$, it follows that $(p_0,q_0)\in C$.
\end{proof}

\begin{proposition}\label{prop:finite}
Let $A$ be a finite poset and let $L$ be a meet-continuous dcpo. Then $[A\to L]$ is meet-continuous.
\end{proposition}

\begin{proof}
We proceed by induction on $|A|$. If $A=\varnothing$, then $[A\to L]$ is a singleton and the conclusion holds.

Suppose that $A$ is nonempty, and choose a maximal element $m\in A$. Let
$$
A_0=A\setminus\{m\},
\qquad
B=\mathord\downarrow_A m\setminus\{m\}.
$$
Then $B$ is a lower subset of $A_0$. By the induction hypothesis, both $[A_0\to L]$ and $[B\to L]$ are meet-continuous.

Define a relation $R\subseteq[A_0\to L]\times L$ by
$$
uRv
\quad\Longleftrightarrow\quad
u(a)\leq v\text{ for every }a\in B.
$$
The relation is order-compatible: if $u'\leq u$, $uRv$ and $v\leq v'$, then
$u'(a)\leq u(a)\leq v\leq v' $ for all $a\in B,$
so $u'Rv'$. It is also nonempty, since for any $c\in L$ the constant map with value $c$ on $A_0$ is related to $c$.

Because $A$ is finite, every directed subset of $A$ has a greatest element. Hence every monotone map from $A$ to $L$ is Scott continuous. A map $f: A\to L$ is monotone precisely when its restriction $u=f|_{A_0}$ is monotone and
$u(a)\leq f(m)$ for all $a\in B.
$
Thus
$
f\longmapsto(f|_{A_0},f(m))
$
is an order isomorphism from $[A\to L]$ onto $R$.

Directed suprema in $R$ exist and are computed coordinatewise. Indeed, if $(u_i,v_i)_{i\in I}$ is a directed family in $R$, then, for each $a\in B$,
$
\bigl(\bigvee_{i\in I} u_i\bigr)(a)
=\bigvee_{i\in I} u_i(a)
\leq\bigvee_{i\in I} v_i.
$
Consequently, $(\bigvee_{i\in I} u_i,\bigvee_{i\in I} v_i)\in R$.

We verify the remaining hypothesis of Lemma~\ref{lem}. Let $D\subseteq L$ be directed, and suppose that $uR\bigvee D$.  For each $d\in D$, let $\kappa_d: B\to L$ be the constant map with value $d$. Then
$u|_B\leq\bigvee_{d\in D}\kappa_d.$
By the meet-continuity of $[B\to L]$,
$$
u|_B\in\cl_{\sigma([B\to L])}W,
\qquad
W=\{h\in[B\to L]: h\leq u|_B,\
h\leq\kappa_d\text{ for some }d\in D\}.
$$

Define $\rho_{u|_B}:[B\to L]\to\mathord\downarrow u|_B$ by
$$
\rho_{u|_B}(h)=
\begin{cases}
h,&h\leq u|_B,\\
u|_B,&h\nleq u|_B.
\end{cases}
$$
The map $\rho_{u|_B}$ is monotone. Indeed, if $h\leq h'$ and $h'\leq u|_B$, then $\rho_{u|_B}(h)=h\leq h'=\rho_{u|_B}(h')$. If $h'\nleq u|_B$, then $\rho_{u|_B}(h')=u|_B$ and $\rho_{u|_B}(h)\leq u|_B$.

It also preserves directed suprema. Let $H\subseteq[B\to L]$ be directed. If every member of $H$ lies below $u|_B$, then $\bigvee H\leq u|_B$, and $\rho_{u|_B}$ fixes every member of $H$ as well as $\bigvee H$. Otherwise, $H$ contains some $h\nleq u|_B$. In this case $\rho_{u|_B}(h)=u|_B$ and $\bigvee H\nleq u|_B$, so both $\bigvee\rho_{u|_B}[H]$ and $\rho_{u|_B}(\bigvee H)$ are equal to $u|_B$. Therefore $\rho_{u|_B}$ is Scott continuous.

For $h\leq u|_B$, define $\psi_u(h):A_0\to L$ by
$$
\psi_u(h)(a)=
\begin{cases}
h(a),&a\in B,\\
u(a),&a\in A_0\setminus B.
\end{cases}
$$
We check that $\psi_u(h)$ is monotone. Since $B$ is lower in $A_0$, if $a\leq a'$ and $a'\in B$, then $a\in B$. Thus the only case in which $a$ and $a'$ lie in different parts is $a\in B$ and $a'\notin B$. In that case,
$$
\psi_u(h)(a)=h(a)\leq u(a)\leq u(a')=\psi_u(h)(a').
$$
Within each part, monotonicity follows from that of $h$ or $u$. Hence $\psi_u(h)$ is monotone and, because $A_0$ is finite, Scott continuous.

The map
$$
\psi_u:\mathord\downarrow u|_B\longrightarrow[A_0\to L]
$$
is Scott continuous: it is monotone, and it preserves directed suprema coordinatewise. Therefore
$$
T=\psi_u\circ\rho_{u|_B}:[B\to L]\longrightarrow[A_0\to L]
$$
is Scott continuous. Also, $T(u|_B)=u$.

Let $h\in W$, and choose $d\in D$ such that $h\leq u|_B,\kappa_d$. Since $h\leq u|_B$, we have $T(h)=\psi_u(h)\leq u$. Moreover,
$
T(h)(a)=h(a)\leq d (a\in B),
$
so $T(h)Rd$. Thus
$$
T[W]\subseteq
\{ b\in[A_0\to L]: b\leq u\text{ and }bRd\text{ for some }d\in D\}.
$$
Applying the Scott-continuous map $T$ to the closure relation for $u|_B$ gives
$$
u=T(u|_B)\in
\cl_{\sigma([A_0\to L])}
\{b\in[A_0\to L]: b\leq u\text{ and }bRd\text{ for some }d\in D\}.
$$
This verifies the hypothesis of Lemma~\ref{lem}. Hence $R$ is meet-continuous. Since $[A\to L]$ is order isomorphic to $R$, it is also meet-continuous.
\end{proof}

We next transfer meet-continuity along approximations of the identity. This gives the bifinite case; Scott-continuous retractions then give the RB case.

\begin{lemma}\label{lem:lifting}
Let $S$ be a dcpo. Let $I$ be a nonempty index set, and suppose that, for every $\alpha \in I$, there are a dcpo $T_\alpha$
and Scott-continuous maps
$S\xrightarrow{p_\alpha}T_\alpha\xrightarrow{j_\alpha}S$ with
$p_\alpha j_\alpha=\id_{T_\alpha}$,
such that 
\begin{enumerate}
    \item $q_\alpha=j_\alpha p_\alpha\leq\id_S$,
    \item the family $(q_\alpha)_\alpha$ is directed,
    \item $\bigvee_\alpha q_\alpha=\id_S$.
\end{enumerate}
If every $T_\alpha$ is meet-continuous, then $S$ is meet-continuous.
\end{lemma}

\begin{proof}
Let $D\subseteq S$ be directed, let $x\leq\bigvee D$, and let $U$ be Scott
open with $x\in U$. Since $x=\bigvee_{\alpha \in I} q_\alpha(x)$, choose $\alpha$
with $q_\alpha(x)\in U$. The Scott-open set $j_\alpha^{-1}(U)$ contains
$p_\alpha(x)$. Meet-continuity of $T_\alpha$, applied to
$p_\alpha(x)\leq\bigvee p_\alpha[D]$, supplies $d\in D$ and
$t\in T_\alpha$ such that $t\leq p_\alpha(x),p_\alpha(d)$ and
$j_\alpha(t)\in U$.
Since $q_\alpha\leq\id_S$, we have $j_\alpha(t)\leq x,d$. Thus
$U\cap\mathord{\downarrow}x\cap\mathord{\downarrow}D\neq\varnothing$,
proving meet-continuity.
\end{proof}

\begin{theorem}\label{thm:rb}
Let $X$ be a bifinite domain, or more generally an RB-domain. If $L$ is a meet-continuous dcpo, then $[X\to L]$ is meet-continuous.
\end{theorem}

\begin{proof}
First suppose that $X$ is bifinite, and let $(h_i)_{i\in I}$ be a directed
approximating identity of idempotent deflations. Let $X_i=h_i[X]$,
let $e_i:X_i\hookrightarrow X$ be inclusion, and let $r_i:X\to X_i$ be the
corestriction of $h_i$. Define
\[
\begin{aligned}
  P_i:[X\to L]&\longrightarrow[X_i\to L],& P_i(f)&=f\circ e_i,\\
  J_i:[X_i\to L]&\longrightarrow[X\to L],& J_i(g)&=g\circ r_i.
\end{aligned}
\]
The image $X_i$ is a finite dcpo. Indeed, if $D\subseteq X_i$ is directed, then Scott continuity and idempotence of $h_i$ give
\[ h_i(\bigvee_XD)=\bigvee_Xh_i[D]=\bigvee_XD,\]
so the supremum belongs to $X_i$. The maps $P_i$ and $J_i$ are Scott-continuous because directed suprema in function spaces are computed pointwise. Then $P_i\circ J_i=\id$, $J_i\circ P_i(f)=f\circ h_i\leq f$, and
$\bigvee_{i\in I}(f\circ h_i)=f$.
The family $(J_i\circ P_i)_{i\in I}$ is directed. Since $X_i$ is finite,
Proposition~\ref{prop:finite} shows that $[X_i\to L]$ is meet-continuous.
Lemma~\ref{lem:lifting} now applies.

Now assume that $X$ is an RB-domain. Since every RB-domain is a Scott-continuous retract of a bifinite domain, we can choose a bifinite domain $Y$ and Scott-continuous maps $X\xrightarrow{e}Y\xrightarrow{r}X$ with $r\circ e=\id_X$. Precomposition gives
$J(f)=f\circ r$ and $R(g)=g\circ e$, with $R\circ J=\id_{[X\to L]}$. Thus $[X\to L]$ is a Scott-continuous retract of the meet-continuous dcpo $[Y\to L]$.
\end{proof}

\begin{example}\label{ex:sequences}
Let $\omega_\infty=\mathbb N\cup\{\infty\}$, ordered by
$0<1<\cdots<\infty$. Since $\omega_\infty$ is an RB-domain, the preceding
theorem implies that $[\omega_\infty\to L]$ is meet-continuous for every
meet-continuous dcpo $L$.
\end{example}

We use the following topologies on a function space \cite{Isbell,XiYang}. 
\begin{definition}
The \emph{pointwise topology} on $[X\to L]$ is the subspace topology inherited from the product $L^X$, where each copy of $L$ is equipped
with its Scott topology. Its subbasic open sets have the form $\{f:f(x)\in U\}$, with $x\in X$ and $U\in\sigma(L)$.
\end{definition}
\begin{definition}
The \emph{Isbell topology} on $[X\to L]$ is generated by the sets $N(\mathcal H,U)=\{f\in[X\to L]:f^{-1}(U)\in\mathcal H\}$, where $U\in\sigma(L)$ and $\mathcal H$ is Scott open in $\mathcal O(X)$. We denote it by $\Is([X\to L])$.
\end{definition}

In general, $\Is([X\to L])\subseteq\sigma([X\to L])$. For a finite set $A$, the Isbell topology on $[A\to L]$ is exactly the pointwise topology. The next example shows that the Isbell topology and the Scott topology on a meet-continuous function space need not coincide.

\begin{example}\label{thm:isbell-scott}
Let $A=\{0,1\}$ be a two-point discrete poset and let $B=\RO(\mathbb R)$ be the complete Boolean algebra of regular open subsets of the real line, ordered by inclusion. Then $B$ is a meet-continuous dcpo, $\Is([A\to B])\subsetneq \sigma([A\to B])$, and $[A\to B]$ is meet-continuous.
\end{example}

In fact, in a complete Boolean algebra, meet with a fixed element is a left adjoint and preserves arbitrary joins; hence $B$ is a meet-continuous dcpo. Since $A$ is discrete, $[A\to B]\cong B\times B$. Because $A$ is finite, the Isbell topology is pointwise, and therefore
$$
  \Is([A\to B])=\sigma(B)\times\sigma(B).
$$

Ern\'e and Gatzke proved that the binary join on $B=\RO(\mathbb R)$ is not continuous when $B\times B$ carries the product of the Scott topologies
\cite{ErneGatzke}. On the other hand, binary join is Scott-continuous on the product dcpo $B\times B$. Equality of the product topology with $\sigma(B\times B)$ would therefore contradict their result. Hence
\[
  \sigma(B)\times\sigma(B)\subsetneq\sigma(B\times B).
\]
%This implication and the role of the regular-open algebra are also recalled in \cite[Section~1]{XiShenZhao}. 
Finally, Proposition~\ref{prop:finite} shows that $[A\to B]$ is meet-continuous.

\section{FS is insufficient}\label{sec:fs}

Theorem~\ref{thm:rb} shows RB-domains are sufficient for the meet-continuity of function spaces for every meet-continuous dcpo. This section tells us that the RB-domain cannot be replaced by the FS-domain. 

Shen, Ho, Xi, and Zhao~\cite[Theorem~4.20]{ShenEtAl} constructed an $\omega$-algebraic domain $L$ for which the function space from $\mathsf{Disk}^{\top}$ to $L$ is not continuous. Using their construction, we prove the failure of the meet-continuity from the closed-disk domain or from domains associated with non-simplicial proper cones.

Let $V$ be a finite-dimensional real vector space. A cone $C\subseteq V$ is called \emph{proper} if it is closed, convex, pointed, and generating. It is \emph{simplicial} if it is generated by a basis of $V$, i.e there is a basis ${b_1, b_2, ..., b_d}$ such that \[
C=\operatorname{cone}\{b_1,\ldots,b_d\}
=\left\{\sum_{i=1}^{d}r_i b_i
\;:dle|\; r_i\ge0,\ i=1,\ldots,d\right\}.
\]
Here pointed means $C\cap(-C)=\{0\}$, and generating means $C-C=V$. The cone order is defined by $x\leq_C y$ if $y-x\in C$. For $x,y\in-C$, its way-below relation is given by
\begin{equation}\label{eq:cone-waybelow}
x\ll y\quad\Longleftrightarrow\quad y-x\in\Int C.
\end{equation}
If $C\subseteq V$ is a proper cone, then
$\Int C\neq\varnothing$ and
$C+\Int C=\Int C$. Let $A\subseteq V$ and let $Q: A\rightarrow V$, $Q$ is \emph{C-monotone} if it preserves the cone order: 
\[
x,y\in A,\quad x\le_C y
\quad\Longrightarrow\quad
Q(x)\le_C Q(y).
\]

The cone domain associated with $C$ is $D_C=(-C)\cup\{\bot\}$, where $-C$ carries the cone order and $\bot$ is a new least element. Keimel \cite[Proposition~5.4]{Keimel2009} showed that $D_C$ is an FS-domain. Chen \cite[Theorem~1.1]{ChenCone} proved that $D_C$ is an RB-domain if and only if $C$ is simplicial.

We shall use the following result of Chen. This is Chen's Corollary~4.8 after choosing linear coordinates on $V$
\cite[Corollary~4.8]{ChenCone}.

\begin{proposition}[{\cite[Corollary~4.8]{ChenCone}}]
\label{prop:cone-obstruction}
Let $C\subseteq V$ be a non-simplicial proper cone and let $\Omega\subseteq V$ be a nonempty open set. There is no sequence of finite-valued $C$-monotone maps $Q_n:\Omega\to V$ converging to the identity uniformly on compact subsets of $\Omega$.
\end{proposition}

For a dcpo $Y$, let
$$
  \mathcal A_{\mathrm{fin}}(Y)
  =\{q\in[Y\to Y]:q[Y]\text{ is finite and }q\leq\id_Y\}.
$$

Let $D_C^\top$ be obtained by adjoining a new greatest element. A directed subset of $D_C^\top$ with supremum $\top$ must contain $\top$, and hence
$\top\ll\top$. It follows that $D_C^\top$ is a countably based continuous domain
with a compact greatest element.

Adjoining this new greatest element also preserves the FS property.
Indeed, extend each map in an FS approximating family by fixing
the new greatest element, and add that element to its finite
separating set. The extensions are Scott-continuous, form a directed
family, and have supremum the identity. Consequently both
$D_C^\top$ and $\mathsf{Disk}^\top$ are countably based FS-domains.

\begin{lemma}
\label{lem:cone-closure}
If $C$ is a non-simplicial proper cone, then \[
\id_{D_C^\top}\notin
  \cl_{\sigma([D_C^\top\to D_C^\top])}\mathcal A_{\mathrm{fin}}(D_C^\top).
\]
\end{lemma}

\begin{proof}
Suppose the contrary. 
Since $-\Int C$ is a nonempty open subset of the finite-dimensional space $V$, it contains a compact full-dimensional parallelotope. Choose a compact full-dimensional parallelotope
$K\subseteq-\Int C$ and an element $a\in\Int C$. Fix $\epsilon>0$. For $x\in K$, let $u_x=x-\frac{3\epsilon}{4}a$, $v_x=x-\frac{\epsilon}{2}a$. Then $u_x, v_x\in-\Int C$ since $C+\Int C\subseteq\Int C$.

By \eqref{eq:cone-waybelow}, $u_x\ll v_x$. Define \[
N_x = \bigl\{ y \in V \mathrel{:}
  x-y+\tfrac{\epsilon}{4}a \in \Int C,\;
  y-x+\tfrac{\epsilon}{2}a \in \Int C
\bigr\}.
\]
Then $K\subseteq \bigcup_{x\in K}N_x$. By compactness of $K$, we can choose $x_1,\ldots,x_m$ such that $K\subseteq\bigcup_{j=1}^mN_{x_j}$. Define $\mathord\Uparrow u=\{z\in D_C^\top: u\ll z\}$, $\mathord\Uparrow u$ is Scott open since $D_C^\top$ is continuous. Since evaluation $\operatorname{eV}$ is continuous, $\operatorname{eV}_{v_{x_j}}^{-1}\bigl(\mathord\Uparrow u_{x_j}\bigr)
=\left\{q:q(v_{x_j})\in\mathord\Uparrow u_{x_j}\right\}$ is Scott open, and
\[
  \mathcal U_\epsilon
  =\bigcap_{j=1}^m\operatorname{eV}_{v_{x_j}}^{-1}\bigl(\mathord\Uparrow u_{x_j}\bigr)
  =\bigcap_{j=1}^m\{q\in[D_C^\top\to D_C^\top]:q(v_{x_j})\in\mathord{\Uparrow}u_{x_j}\}
\]
is a Scott-open neighborhood of $\id_{D_C^\top}$. Since $\id_{D_C^\top}\in
  \cl_{\sigma([D_C^\top\to D_C^\top])}\mathcal A_{\mathrm{fin}}(D_C^\top)$, $\mathcal U_\epsilon\cap\mathcal A_{\mathrm{fin}}(D_C^\top)\neq\varnothing$ . Hence there is
$q_\epsilon\in\mathcal{U_\epsilon}\cap
\mathcal A_{\mathrm{fin}}(D_C^\top)$.

If $y\in K\cap N_{x_j}$, then $y-\epsilon a\le_C u_{x_j} $ and $v_{x_j}\le_C y$. Since $q_\epsilon\in\mathcal{U_\epsilon}$, $u_{x_j}\le q_\epsilon(v_{x_j})$. The monotonicity of $q_\epsilon$ and $q_\epsilon\leq\id_{D_C^\top}$ give
\begin{equation}
  y-\epsilon a
  \leq_C u_{x_j}
  \leq_C q_\epsilon(v_{x_j})
  \leq_C q_\epsilon(y)
  \leq_C y.
  \label{eq:cone-Q-bound}
\end{equation}
In particular, $q_\epsilon(y)$ is a non-bottom element of $D_C$.
Therefore $Q_\epsilon=q_\epsilon|_K$ is a finite-valued $C$-monotone map into $V$. By \eqref{eq:cone-Q-bound}, $y-Q_\epsilon(y)\in[0,\epsilon a]_C=\epsilon[0,a]_C$.
The order interval $[0,a]_C=\{z: z-0\in C, a-z\in C\}=C\cap (a-C)$ is compact, so
\[
  \sup_{y\in K}\lVert Q_\epsilon(y)-y\rVert
  \leq M_a\epsilon
\]
for a constant $M_a=\max_{z\in[0,a]_C}\lVert z\rVert<\infty$ independent of $\epsilon$. Taking
$\epsilon_n\searrow0$ and setting
$Q_n=Q_{\epsilon_n}|_{\Int K}$ gives finite-valued $C$-monotone maps
$Q_n:\Int K\to V$ converging uniformly on compact subsets of $\Int K$ to the
identity. This contradicts Proposition~\ref{prop:cone-obstruction}.
\end{proof}

Recall that $\mathsf{Disk}$ consists of all closed disks $\overline B(z,r)\subseteq\mathbb R^2$ with $r\geq0$, together with the whole plane as its least element, ordered by reverse inclusion. Let
\[
  C=\{(t,z)\in\mathbb R\times\mathbb R^2:t\geq\lVert z\rVert\}
\]
be the three-dimensional Lorentz cone, and let $H=(-\infty,0]\times\mathbb R^2$. The map $\overline B(z,r)\mapsto(-r,z)$ is a bijection from the non-bottom part of $\mathsf{Disk}$ onto $H$. It preserves and reflects the order, since
\[
\begin{aligned}
  \overline B(z,r)\leq_{\mathsf{Disk}}\overline B(w,s)
  &\Longleftrightarrow
    \overline B(w,s)\subseteq\overline B(z,r)\\
  &\Longleftrightarrow \lVert z-w\rVert+s\leq r\\
  &\Longleftrightarrow (r-s,w-z)\in C\\
  &\Longleftrightarrow (-r,z)\leq_C(-s,w).
\end{aligned}
\]
Thus $\mathsf{Disk}\cong(H,\leq_C)_\bot$ \cite[Section~3]{ChenKouLyu}, and the way-below relation between non-bottom elements is given by \eqref{eq:cone-waybelow}.

Let $p_0=\overline B(0,0)$ and $I=\downarrow_{\mathsf{Disk}}p_0$. For a closed disk $\overline B(z,r)$, we have
\[
\begin{aligned}
  \overline B(z,r)\in I
  &\Longleftrightarrow 0\in\overline B(z,r)\\
  &\Longleftrightarrow r\geq\lVert z\rVert
  \Longleftrightarrow (-r,z)\in-C.
\end{aligned}
\]
Hence the disk representation restricts to an order isomorphism between $I\setminus\{\bot\}$ and $-C$. Including the least element gives $I\cong D_C$.

The Lorentz cone $C$ is proper \cite[Section~3]{ChenKouLyu} and non-simplicial \cite[Example~5.4]{ChenCone}. The domain $\mathsf{Disk}$ is an FS-domain \cite[Section~4.2.2]{AJ} (see also \cite{Lawson2008}). It is not an RB-domain \cite[Theorem~5.6]{ChenKouLyu}. After adjoining a new greatest element, $\mathsf{Disk}^\top$ remains an FS-domain, and its new greatest element is compact.

\begin{corollary}
\label{cor:disk-closure}
Let $\mathsf{Disk}^\top=\mathsf{Disk}\cup\{\top\}$. Then 
\[
  \id_{\mathsf{Disk}^\top}\notin
  \cl_{\sigma([\mathsf{Disk}^\top\to \mathsf{Disk}^\top])}
  \mathcal A_{\mathrm{fin}}(\mathsf{Disk}^\top).
\]
\end{corollary}

\begin{proof}
Identify $D_C^\top$ with $I\cup\{\top\}$ in $\mathsf{Disk}^\top$.
Let $i:D_C^\top\to \mathsf{Disk}^\top$ be the inclusion and define
$r:\mathsf{Disk}^\top\to D_C^\top$ by
\[
  r(x)=
  \begin{cases}
    x, & x\in I,\\
    \top, & x\notin I.
  \end{cases}
\]
For every directed $E\subseteq D_C^\top$, its supremum in
$\mathsf{Disk}^\top$ lies in $I$ if $E\subseteq I$, and equals
$\top$ if $\top\in E$. Thus $i$ preserves directed suprema.
By construction, $r(x)$ is the least element of $D_C^\top$
above $x$, so
\[
  r(x)\leq u\quad\Longleftrightarrow\quad x\leq i(u)
  \qquad(x\in \mathsf{Disk}^\top,\ u\in D_C^\top).
\]
Therefore $r$ is Scott-continuous as a left adjoint, and
$r\circ i=\id_{D_C^\top}$.

Define $\Phi:[\mathsf{Disk}^\top\to \mathsf{Disk}^\top]
\longrightarrow [D_C^\top\to D_C^\top]$, $\Phi(q)=r\circ q\circ i$. The map $\Phi$ is monotone and preserves directed suprema, since
these suprema are computed pointwise and $r$ is Scott-continuous. For $q\in\mathcal A_{\mathrm{fin}}(\mathsf{Disk}^\top)$, we have 
\[
\Phi(q)[D_C^\top]\subseteq r[q[\mathsf{Disk}^\top]], \quad \Phi(q)\leq r\circ i=\id_{D_C^\top}.
\] 
The set $r[q[\mathsf{Disk}^\top]]$ is finite. Consequently, 
\[\Phi(\id_{\mathsf{Disk}^\top})=\id_{D_C^\top}, \quad \Phi[\mathcal A_{\mathrm{fin}}(\mathsf{Disk}^\top)]\subseteq\mathcal A_{\mathrm{fin}}(D_C^\top).
\] 
If $\id_{\mathsf{Disk}^\top}\in\cl_{\sigma([\mathsf{Disk}^\top\to \mathsf{Disk}^\top])} \mathcal A_{\mathrm{fin}}(\mathsf{Disk}^\top)$, Scott continuity of $\Phi$ would therefore give $\id_{D_C^\top}\in\cl_{\sigma([D_C^\top\to D_C^\top])}\mathcal A_{\mathrm{fin}}(D_C^\top)$, contrary to Lemma~\ref{lem:cone-closure}.
\end{proof}

The following result comes from Shen, Ho, Xi, and Zhao\cite[Proposition~4.9, Lemma~4.10, and Proposition~4.16]{ShenEtAl}. For a continuous domain $Y$, choose a basis $B$ of $Y$ containing its greatest element, and set
\[
P=\{F\subseteq B:F\text{ is nonempty, finite, and directed}\},\qquad L_Y=\Idl(P),
\]
where $P$ is ordered by inclusion and ideals are nonempty directed lower sets. If $Y$ has a greatest element $\top_Y$, for $F,G\in P$,
the set $F\cup G\cup\{\top_Y\}$ belongs to $P$ and contains both. Thus $P$ is directed and is the greatest element of $L_Y$.

\begin{lemma}\label{lem:algebraic-cover}
Let $Y$ be a continuous domain with a greatest element. There is an algebraic domain $L_Y$ with a greatest element and Scott-continuous maps $Y\xrightarrow{s}L_Y\xrightarrow{r}Y$, $r\circ s=\id_Y$,
such that $\downarrow_{L_Y}k$ is finite for every compact $k\in L_Y$.
Consequently,
\begin{equation}\label{eq:finite-downset}
a\ll_{L_Y}\top_{L_Y}\quad\Longrightarrow\quad
\downarrow_{L_Y}a\text{ is finite}.
\end{equation}
If $Y$ is countably based, $L_Y$ may be chosen $\omega$-algebraic. If the greatest element of $Y$ is compact, the construction also satisfies $s(\top_Y)=\top_{L_Y}$.
\end{lemma}

\begin{proposition}\label{prop:closure-obstruction}
Let $Y$ be a continuous domain with a greatest element. If
\[
\id_Y\notin\cl_{\sigma([Y\to Y])}\Afin(Y),
\]
then $[Y\to L_Y]$ is not meet-continuous.
\end{proposition}

\begin{proof}
Suppose that $[Y\to L_Y]$ were meet-continuous. Write $L_Y=\Idl(P)$ as in the construction above. For each $F\in P$, let $c_F\in [Y\to L_Y]$ be the constant map with value $\mathord\downarrow_P F$. The family $\mathcal D=\{c_F:F\in P\}$ is directed, with supremum the constant map of value $\top_{L_Y}$. The maps $r, s$ of Lemma~\ref{lem:algebraic-cover} are defined by $r(I)=\bigvee(\bigcup_{F\in I}F)$, $s(y)=\{F\in P:F\subseteq B\cap\mathord\Downarrow y\}$ where ${\mathord\Downarrow y=\{z\in Y: z\ll y}\}$. Since $s\le\bigvee\mathcal D$, meet-continuity gives
\begin{equation}\label{eq:section-closure}
s\in\cl_{\sigma([Y\to L_Y])}
(\mathord\downarrow s\cap\mathord\downarrow\mathcal D).
\end{equation}
If $u\le s, c_F$, then $u[Y]\subseteq\mathord\downarrow_{L_Y}(\mathord\downarrow_P F)$. Since $F$ is finite, $\mathord\downarrow_{L_Y}(\mathord\downarrow_P F)$ is finite and hence $u$ has finite image. Therefore $r\circ u$ has finite image and $r\circ u\le r\circ s=\id_Y$. Define $R: [Y\to L_Y]\longrightarrow[Y\to Y]$ with $R(f)=r\circ f$. Clearly $R$ is Scott continuous. Applying it to \eqref{eq:section-closure} yields  $\id_Y\in\cl_{\sigma([Y\to Y])}\mathcal A_{\mathrm{fin}}(Y)$, a contradiction.
\end{proof}

\begin{theorem}
\label{thm:fs-with-top}
For every non-simplicial finite-dimensional proper cone $C$, there is an $\omega$-algebraic domain $L_{D_C^\top}$ such that $[D_C^\top\to L_{D_C^\top}]$ is not meet-continuous.
\end{theorem}

\begin{proof}
We first show that $D_C$ has a countable basis. Indeed, $-\Int C\subseteq V$ is a nonempty open set. It is separable metric space and thus there exists a countably dense subset $B\subseteq -\Int C$. For each $y\in -C$, fix $a\in \Int C$. Let $y_n=y-\frac{1}{n}a$, then $y_n\in-\Int C$ since $C+\Int C\subseteq \Int C$, and $y_{n+1}-y_{n}=(\frac{1}{n}-\frac{1}{n+1})a\in C$, $y-y_n=\frac{1}{n}a\in\Int C$. Hence $y_1\le_C y_2\le_C \cdots\le_C y$. If $z\in D_C$ is an upper bound of $(y_n)$, then $z\in-C$ and $z-y_n\in C$ for every $n$.Since $C$ is closed and $y_n\to y$, we have $z-y\in C$. So $y=\bigvee_n y_n$. For each $n\in \mathrm N$, define $W_n=(y_n+\Int C)\cap(y-\Int C)$ and it is nonempty open since $y-\frac{1}{2n}a\in W_n$. Since $W_n\subseteq -\Int C$ and $B$ is dense, there is $b_n\in B\cap W_n$ with $y_n\le_C b_n\ll_C y$. Then 
\[
y\le_C \bigvee\{b_n\in B: y_n\le_C b_n\ll_C y\}\le_C\bigvee\{b\in B: b\ll_C y\}\le_C y.
\]
Conversely $y=\bigvee\{b\in B: b\ll_C y\}$. The set $\{b\in B: b\ll_C y\}\cup\{\bot\}$ is directed clearly by the directness of $(y_n)_n$. Thus $D_C$ has a countable basis. 

Since $C$ is pointed and $D_C$ is an FS-domain, every directed subset of $D_C$ has supremum and hence $\top\in D_C^{\top}$ is a compact element. So $D_C^{\top}$ has a countable basis.
Apply Lemma~\ref{lem:algebraic-cover} to $Y=D_C^\top$ and denote the resulting $\omega$-algebraic domain by $L_{D_C}$.
Lemma~\ref{lem:cone-closure} gives the hypothesis of Proposition~\ref{prop:closure-obstruction}, which implies that $[D_C^\top\to L_{D_C^\top}]$ is not meet-continuous.
\end{proof}

\begin{corollary}
 There is likewise an $\omega$-algebraic domain $L_{\mathsf Disk^\top}$ such that $[\mathsf{Disk}^\top\to L_{\mathsf Disk^\top}]$ is not meet-continuous.
\end{corollary}

\begin{proof}
The disks with rational centres and positive rational radii i.e. $\{\overline B(z,r): z\in \mathrm Q^2, r\in \mathrm Q_+\}$, together with the least element and the new greatest element, form a countable basis of $\mathsf{Disk}^{\top}$. Apply Lemma~\ref{lem:algebraic-cover} to $Y=\mathsf{Disk}^\top$ to obtain $L_{\mathsf{Disk}^\top}$. The conclusion follows from Corollary~\ref{cor:disk-closure} and Proposition~\ref{prop:closure-obstruction}.
\end{proof}

\section{Characterize RB-domains by the meet-continuity}\label{sec:necessary}
Theorem~\ref{thm:rb} shows that, if $X$ is an RB-domain, then $[X\to L]$ is meet-continuous for every meet-continuous dcpo $L$. As the following example shows, for a particular choice of $L$, the function space $[X\to L]$ can be meet-continuous even when $X$ is not an RB-domain. In this section, we complete the characterization by proving that a continuous domain $X$ is an RB-domain whenever $[X\to L]$ is meet-continuous for every meet-continuous dcpo $L$.

\begin{example}
Let $L_0=\{\alpha, \beta, \top\}$ where $\alpha\le \top, \beta\le\top$ and $\alpha\parallel \beta$. Let $B=\RO(\mathbb R)$. Then $B$ is not continuous, so $B$ is not an RB-domain. However $[B\to L_0]$ is meet-continuous since  $B$ is Scott compact. In fact, it can be proved that for any Scott compact dcpo $X$, the function space $[X\to L_0]$ is always meet-continuous.
\end{example}

Let $P_1$ and $P_2$ be dcpos. Their \emph{ordinal sum} $P_1+P_2$ is the disjoint union $P_1\sqcup P_2$, ordered by keeping the original order on each summand and placing every element of $P_1$ strictly below every element of $P_2$.
It is a dcpo:
let $D\subseteq P_1+P_2$ be directed. If $D\cap P_2\neq\varnothing$, then $D\cap P_2$ is cofinal in $D$ and
\[
\bigvee_{P_1+P_2}D=\bigvee_{P_2}(D\cap P_2).
\]
If $D\subseteq P_1$, then
\[
\bigvee_{P_1+P_2}D=\bigvee_{P_1}D.
\]

For $u\in[P_1\to L]$ and $v\in[P_2\to L]$, write $uR v$ if
$u(p)\leq v(q)$ for all $p\in P_1$ and $q\in P_2$. If
$f\in[P_1+P_2\to L]$, then its restrictions
$u=f|_{P_1}$ and $v=f|_{P_2}$ are Scott-continuous. Since every element
of $P_1$ lies below every element of $P_2$, monotonicity of $f$ gives
$u(p)\leq v(q)$ for all $p\in P_1$ and $q\in P_2$; hence $uR v$.

Conversely, suppose that $u\in[P_1\to L]$, $v\in[P_2\to L]$, and
$uR v$. Define $f:P_1+P_2\to L$ by
\[
f(x)=
\begin{cases}
u(x), & x\in P_1,\\
v(x), & x\in P_2.
\end{cases}
\]
The monotonicity of $u$ and $v$, together with $uR v$, shows that
$f$ is monotone. Let $D\subseteq P_1+P_2$ be directed. If $D\subseteq
P_1$, then $\bigvee D$ is computed in $P_1$, and the Scott continuity of
$u$ gives $f(\bigvee D)=\bigvee f(D)$. If $D\cap P_2\neq\varnothing$,
then $D\cap P_2$ is cofinal in $D$ and
$\bigvee D=\bigvee_{P_2}(D\cap P_2)$. Therefore, by the Scott continuity
of $v$,
\[
f(\bigvee D)
=v\left(\bigvee(D\cap P_2)\right)
=\bigvee v(D\cap P_2)
=\bigvee f(D).
\]
Thus $f$ is Scott-continuous. The restriction map and the preceding
piecewise construction are mutually inverse and preserve the pointwise
orders. Hence
\begin{equation}\label{eq:ordinal-iso}
[P_1+P_2\to L]
\cong
\{(u,v)\in[P_1\to L]\times[P_2\to L]:uR v\}.
\end{equation}

\begin{theorem}\label{thm:ordinal}
Let $P_1,P_2,L$ be dcpos, and suppose that $P_2$ has a least element $q_0$. If both $[P_1\to L]$ and
$[P_2\to L]$ are meet-continuous, then $[P_1+P_2\to L]$ is
meet-continuous.
\end{theorem}

\begin{proof}
Under \eqref{eq:ordinal-iso}, the relation $R$ is order-compatible and
the associated sub-dcpo has coordinatewise directed suprema. Let
$D\subseteq[P_2\to L]$ be directed, and suppose that
$v\leq\bigvee D$ and $uR v$. For $d\in D$, let
$\kappa_d:P_1\to L$ be constant with value $d(q_0)$. If $b\leq\kappa_d$,
then $b(p)\leq d(q_0)\leq d(q)$ for all $p\in P_1$ and $q\in P_2$, so
$bR d$. Moreover,
$u(p)\leq v(q_0)\leq\bigvee_{d\in D}d(q_0)=
\bigvee_{d\in D}\kappa_d(p)$.
Thus $u\leq\bigvee_{d\in D}\kappa_d$. Meet-continuity of $[P_1\to L]$
gives
\[
  u\in\cl_{\sigma([P_1\to L])}
  \bigl(\mathord{\downarrow}u
  \cap\mathord{\downarrow}\{\kappa_d:d\in D\}\bigr).
\]
Every map in the set inside the closure is below $u$ and is related by
$R$ to some $d\in D$. Hence the additional closure hypothesis of
Lemma~\ref{lem} holds, and that lemma applies.
\end{proof}

\begin{remark}\label{rem:removing-top}
For a dcpo $X$, the notation $X+1$ denotes the ordinal sum of $X$ and a one-element dcpo. Thus $X+1$ is obtained from $X$ by adjoining a new greatest element.

Let $L$ be a meet-continuous dcpo. Since $[1\to L]\cong L$, the ordinal-sum theorem gives
\[
[X\to L]\text{ is meet-continuous}
\quad\Longrightarrow\quad
[X+1\to L]\text{ is meet-continuous}.
\]
Consequently, if $[X+1\to L]$ is not meet-continuous, then neither is $[X\to L]$.

In particular, for the cone domains considered above, $D_C=(-C)\cup\{\bot\}$ and $D_C^\top=D_C+1$. Likewise, $\mathsf{Disk}^\top=\mathsf{Disk}+1$. Hence $[D_C\to L_{D_C}]$ and $[\mathsf{Disk}\to L_{\mathsf{Disk}}]$ both are not meet-continuous.
\end{remark}

\begin{lemma}\label{lem:deflation-closure}
Let $Y$ be a continuous domain. If $\id_Y\in\cl_{\sigma([Y\to Y])}\Afin(Y)$, where
  $\mathcal A_{\mathrm{fin}}(Y)
  =\{q\in[Y\to Y]:q[Y]\text{ is finite and }q\leq\id_Y\},$ then $Y$ is an RB-domain.
\end{lemma}

\begin{proof}
It suffices to look for a directed family of deflations with $\id_Y$ as its supremum. Consider $D=\{q^2 : q\in\Afin(Y)\}$, then $D\neq\varnothing$, since $\id_Y\in\cl_{\sigma([Y\to Y])}\Afin(Y)$, $\Afin(Y)\neq\varnothing$. The members of $D$ have finite image and lie below $\id_Y$, hence are deflations. 

\vspace{0.2cm}
\noindent{\bf Claim 1.} $D$ is directed.

If $q\in\Afin(Y)$ and $y\in Y$, then $q(y)\ll y$. Indeed, $q(y)=\bigvee_{z\ll y}q(z)$. The set $\{q(z): z\ll y\}$ is directed and finite, so $q(y)=q(z)$ for some $z\ll y$. Thus $q(y)=q(z)\le z\ll y$.

For $q\in\Afin(Y)$, define
\[
U_q=\bigcap_{s\in q[Y]}\{h\in[Y\to Y]:q(s)\ll h(s)\}.
\]
By the similarly prove of \ref{lem:cone-closure}, $U_q=\bigcap_{s\in q[Y]}\operatorname{eV}_s^{-1}(\mathcal\Uparrow q(s))$. Continuity of $Y$, Scott continuity of evaluation, and finiteness of $q[Y]$ show that $U_q$ is Scott open. The preceding observation gives $\id_Y\in U_q$. Moreover, if $h\in U_q$, then
\[
q^2(y)=q(q(y))\le h(q(y))\le h(y)\qquad(y\in Y),
\]
so $q^2\le h$.
The map $S:[Y\to Y]\to[Y\to Y]$, defined by $S(h)=h^2$, is Scott continuous. It is monotone, and for every directed family $(h_i)_{i\in I}$,
\[
\left(\bigvee_i h_i\right)^2=\bigvee_{i,j}h_i\circ h_j=\bigvee_i h_i^2.
\]
The last equality follows by choosing, for each $i,j$, an index above both.

Let $q,r\in\Afin(Y)$. The set $S^{-1}(U_q\cap U_r)$ is a Scott-open neighbourhood of $\id_Y$. The hypothesis $\id_Y\in\cl_{\sigma([Y\to Y])}\Afin(Y)$ implies $S^{-1}(U_q\cap U_r)\cap \Afin(Y)\neq\varnothing$. Hence there is $t\in\Afin(Y)$ with $t^2\in U_q\cap U_r$, and so $q^2,r^2\le t^2$. Thus $D=\{q^2:q\in\Afin(Y)\}$ is directed. 

\vspace{0.2cm}

\noindent{\bf Claim 2.} $\bigvee D=\id_Y$.

Scott continuity of $S$ gives 
\[
\id_Y=S(\id_Y)\in \cl_{\sigma([Y\to Y])}S(\Afin(Y))=\cl_{\sigma([Y\to Y])}D.
\] 
Since $\downarrow\bigvee D$ is Scott closed and every member of $D$ is below $\id_Y$, it follows that $\id_Y\le\bigvee D\le\id_Y$. 

Thus $Y$ is an RB-domain.
\end{proof}

\begin{theorem}\label{thm:necessary}
Let $X$ be a continuous domain. If $[X\to A]$ is meet-continuous for every
algebraic domain $A$, then $X$ is an RB-domain.
\end{theorem}

\begin{proof}
Adjoin a new greatest element to $X$ and write $Y=X+1$. Then $Y$ is continuous. Let $A=L_Y$ be the algebraic domain supplied by Lemma~\ref{lem:algebraic-cover} for $Y$.

By hypothesis, $[X\to L_Y]$ is meet-continuous.
Since $L_Y$ is algebraic, $[1\to L_Y]\cong L_Y$ is also meet-continuous. Theorem~\ref{thm:ordinal} therefore
shows that $[Y\to L_Y]$ is meet-continuous. Proposition~\ref{prop:closure-obstruction} gives
\[
\id_Y\in
\cl_{\sigma([Y\to Y])}\Afin(Y).
\]
By Lemma~\ref{lem:deflation-closure}, $Y$ is an RB-domain.

Choose a directed family $(q_i)_{i\in I}$ of deflations
on $Y$ with supremum $\id_Y$.
For $x\in X$, the inequality $q_i(x)\leq x$ implies
$q_i(x)\in X$. Directed suprema of subsets of $X$
are unchanged in $Y$, so each restriction
$q_i|_X:X\to X$ is Scott-continuous.
These restrictions form a directed family of deflations
with supremum $\id_X$. Hence $X$ is an RB-domain.
\end{proof}

\begin{theorem}\label{thm:characterization}
For a continuous domain $X$, the following conditions are equivalent.
\begin{enumerate}
\renewcommand{\labelenumi}{(\roman{enumi})}
\item $X$ is an RB-domain.
\item $[X\to L]$ is meet-continuous for every meet-continuous dcpo $L$.
\item $[X\to A]$ is meet-continuous for every algebraic domain $A$.
\end{enumerate}
\end{theorem}

\begin{proof}
The implication (i)$\Rightarrow$(ii) is Theorem~\ref{thm:rb}. Every algebraic domain is meet-continuous, so
(ii)$\Rightarrow$(iii). Theorem~\ref{thm:necessary} proves (iii)$\Rightarrow$(i).
\end{proof}

\begin{remark}
Chen characterizes RB-domains by continuity of $[X\to E]$ for every continuous domain $E$\cite{ChenFunction}. Theorem~\ref{thm:characterization} concerns meet-continuity of the function spaces. In particular, its third condition requires only meet-continuity for algebraic domains.
\end{remark}

For algebraic domains, the characterization can be expressed using bifiniteness. 

\begin{corollary}\label{cor:algebraic}
For an algebraic domain $X$, the following conditions are equivalent.
\begin{enumerate}
\renewcommand{\labelenumi}{(\roman{enumi})}
\item $X$ is bifinite.
\item $X$ is an RB-domain.
\item $[X\to L]$ is meet-continuous for every meet-continuous dcpo $L$.
\item $[X\to A]$ is meet-continuous for every algebraic domain $A$.
\end{enumerate}
\end{corollary}

\begin{proof}
Since $X$ is algebraic, it is continuous.
Theorem~\ref{thm:characterization} therefore gives the equivalence of (ii), (iii), and (iv). Every bifinite domain is an RB-domain by definition. Conversely, every algebraic RB-domain is bifinite in the sense used here \cite[Corollary~4.3 and Theorem~3.11]{JungCCC}. Thus all four conditions are equivalent.
\end{proof}

Finally, we introduce a necessary condition of meet-continuity of function spaces $[X\to L]$ for every meet-continuous dcpo $L$.

\begin{proposition}\label{prop:lawsoncompact}
For continuous domain $X$, if $[X\to L]$ is meet-continuous for every meet-continuous dcpo $L$, then $X$ is Lawson compact.
\end{proposition}

\begin{proof}
Let $L_0=\{\alpha, \beta, \top\}$ where $\alpha\le \top, \beta\le\top$ and $\alpha\parallel \beta$. Let $L_6=\{0,c,d,a,b,1\}$, $0<c,d<a,b<1$, where $c\parallel d$ and $a\parallel b$. 

\begin{center}
\tikzset{
  posetvertex/.style={
    circle,
    fill=black,
    inner sep=1.5pt
  }
}
\begin{tabular}{@{}c@{\hspace{2.5cm}}c@{}}
% L_0
\begin{tikzpicture}[baseline=(current bounding box.center)]
  \node[posetvertex,label=below:{$\alpha$}]
    (alpha) at (-0.8,0) {};
  \node[posetvertex,label=below:{$\beta$}]
    (beta) at (0.8,0) {};
  \node[posetvertex,label=above:{$\top$}]
    (top) at (0,1.4) {};

  \draw (alpha)--(top)--(beta);
\end{tikzpicture}
&
% L_6
\begin{tikzpicture}[baseline=(current bounding box.center)]
  \node[posetvertex,label=below:{$0$}]
    (zero) at (0,0) {};
  \node[posetvertex,label=left:{$c$}]
    (c) at (-0.8,1) {};
  \node[posetvertex,label=right:{$d$}]
    (d) at (0.8,1) {};
  \node[posetvertex,label=left:{$a$}]
    (a) at (-0.8,2) {};
  \node[posetvertex,label=right:{$b$}]
    (b) at (0.8,2) {};
  \node[posetvertex,label=above:{$1$}]
    (one) at (0,3) {};

  \draw (zero)--(c)
        (zero)--(d)
        (c)--(a)
        (c)--(b)
        (d)--(a)
        (d)--(b)
        (a)--(one)
        (b)--(one);
\end{tikzpicture}
\\[0.8em]
$L_0$ & $L_6$
\end{tabular}
\end{center}

\vspace{0.2cm}
\noindent{\bf Claim 1.} 
If $[X\to L_0]$ is meet-continuous, then $X$ is Scott compact. 

Suppose that $X$ is not Scott compact. Then there exists a directed family $(W_i)_{i\in I}$ of Scott-open subsets of $X$ such that $\bigcup_{i\in I}W_i=X$ but every finite subsets of $(W_i)_{i\in I}$ can not cover $X$.
Define
$g_i:X\to L_0$ by
\[
g_i(x)=
\begin{cases}
\top, & x\in W_i,\\
\beta, & x\notin W_i.
\end{cases}
\]
Each $g_i$ is Scott-continuous: the inverse images of $\{\top\}$ and $\{\alpha,\top\}$ are $W_i$, and that of$\{\beta,\top\}$ is $X$. Then $(g_i)_{i\in I}$ is directed and $\bigvee_i g_i=c_\top$, where $c_\top$ is a constant map with $c_\top(x)=\top$ for each $x\in X$. Similarly define $c_\alpha$ is constant map with image $\alpha$. Then $c_\alpha\leq c_\top$, while $c_\alpha\nleq g_i$ for every $i$, because $W_i\neq X$ and $\alpha$ and $\beta$ are incomparable. Since $c_\alpha$ is minimal, $\mathord{\downarrow}c_\alpha\cap\mathord{\downarrow}\{g_i:i\in I\}=\varnothing$. Thus
\[
c_\alpha\leq\bigvee_i g_i\quad\text{but}\quad c_\alpha\notin\cl_{\sigma([X\to L_0])}\left(\mathord{\downarrow}c_\alpha\cap\mathord{\downarrow}\{g_i:i\in I\}\right),
\]
contradicting the meet-continuity. Hence $X$ is Scott compact.

\vspace{0.2cm}
\noindent{\bf Claim 2.} 
If $[X\to L_6]$ is meet-continuous, then $X$ is coherent. 

Suppose that $X$ is not coherent. $X$ is continuous, thus is a locally compact sober dcpo. $L_6$ is a pointed bicomplete sober dcpo but not an L-dcpo since $\mathord\downarrow 1$ is not a complete lattice. Hence $[X\to L_6]$ is not meet-continuous by \cite[Theorem~6.3]{JiaJungLiDichotomy}, a contradiction.

We have shown that $X$ is Scott compact and coherent.
Since $X$ is continuous, it is also well-filtered by
\cite[Proposition~4.2.14(2)]{AJ}. Therefore $X$ is Lawson compact by \cite[Theorem~3.4]{JiaJungLiCoherence}.
\end{proof}

\medskip
\noindent\textbf{Declaration of Generative AI and AI-assisted Technologies}\quad
Generative AI and AI-assisted Technologies were used under the authors' guidance to assist with certain technical steps in the proofs. All principal mathematical ideas and proof strategies were developed by the authors, who independently verified all AI-assisted arguments and take full responsibility for the results.


\begin{thebibliography}{99}
\bibitem{Abbadini}
M.~Abbadini,
\newblock Formal squares over the unit square separate FS-domains from RB-domains,
\newblock arXiv preprint arXiv:2608.08681v1, 2026.

\bibitem{AJ}
S.~Abramsky and A.~Jung,
\newblock Domain theory,
\newblock in S.~Abramsky, D.~M.~Gabbay, and T.~S.~E.~Maibaum (eds.),
\emph{Handbook of Logic in Computer Science}, Vol.~3,
Clarendon Press, Oxford, 1994, pp.~1--168.

\bibitem{ChenFunction}
Y.~Chen,
\newblock A function space characterization of RB-domains,
\newblock arXiv preprint arXiv:2609.20187v1, 2026.

\bibitem{ChenCone}
Y.~Chen,
\newblock Cone domains separate FS-domains from RB-domains,
\newblock arXiv preprint arXiv:2607.02251v1, 2026.

\bibitem{ChenKouLyu}
Y.~Chen, H.~Kou, and Z.~Lyu,
\newblock FS-domains are not always RB-domains,
\newblock arXiv preprint arXiv:2607.00568v2, 2026.

\bibitem{ErneGatzke}
M.~Ern\'e and H.~Gatzke,
\newblock Convergence and continuity in partially ordered sets and semilattices,
\newblock in R.-E.~Hoffmann and K.~H.~Hofmann (eds.),
\emph{Continuous Lattices and Their Applications},
Lecture Notes in Pure and Applied Mathematics, Vol.~101, Marcel Dekker, New York, 1985, pp.~9--40.

\bibitem{Gierz}
G.~Gierz, K.~H.~Hofmann, K.~Keimel, J.~D.~Lawson, M.~Mislove, and D.~S.~Scott,
\newblock \emph{Continuous Lattices and Domains},
\newblock Encyclopedia of Mathematics and its Applications, Vol.~93, Cambridge University Press, Cambridge, 2003.

\bibitem{Isbell}
J.~R.~Isbell,
\newblock Function spaces and adjoints,
\newblock \emph{Mathematica Scandinavica}
36 (1975), 317--339.

\bibitem{JiaThesis}
X.~Jia,
\newblock \emph{Meet-Continuity and Locally Compact Sober Dcpos},
\newblock Ph.D. thesis, University of Birmingham, 2018.


\bibitem{JiaJungLiCoherence}
X.~Jia, A.~Jung, and Q.~Li,
\newblock A note on coherence of dcpos,
\newblock \emph{Topology and its Applications}
209 (2016), 235--238.

\bibitem{JiaJungLiDichotomy}
X.~Jia, A.~Jung, and Q.~Li,
\newblock A dichotomy result for locally compact sober dcpos,
\newblock \emph{Houston Journal of Mathematics}
45 (2019), no.~3, 935--951.

\bibitem{JungCCC}
A.~Jung,
\newblock \emph{Cartesian Closed Categories of Domains},
\newblock CWI Tracts, Vol.~66,
Centrum voor Wiskunde en Informatica, Amsterdam, 1989.

\bibitem{Keimel2009}
K.~Keimel,
\newblock Bicontinuous domains and some old problems in domain theory,
\newblock \emph{Electronic Notes in Theoretical Computer Science}
257 (2009), 35--54.

\bibitem{KouLiuLuo}
H.~Kou, Y.-M.~Liu, and M.-K.~Luo,
\newblock On meet-continuous dcpos,
\newblock in G.~Q.~Zhang, J.~D.~Lawson, Y.-M.~Liu,
and M.-K.~Luo (eds.),
\emph{Domain Theory, Logic and Computation},
Semantic Structures in Computation, Vol.~3,
Springer, Dordrecht, 2003, pp.~117--135.

\bibitem{Lawson2008}
J.~D.~Lawson,
\newblock Metric spaces and FS-domains,
\newblock \emph{Theoretical Computer Science}
405 (2008), no.~1--2, 73--74.

\bibitem{LiuEtAl2020}
D.~Liu, X.~Zhou, C.~Lu, and Q.~Li,
\newblock Coincidence of the Isbell and Scott topologies on the function spaces of quasicontinuous domains,
\newblock \emph{Topology and its Applications}
285 (2020), 107407.

\bibitem{LuanLi2025}
W.~Luan and Q.~Li,
\newblock Function spaces of Lawson-compact algebraic domains,
\newblock \emph{Topology and its Applications}
366 (2025), 109295.

\bibitem{MiaoEtAl}
H.~Miao, W.~Luan, D.~Liu, and Q.~Li,
\newblock Spaces of functions from RB-domains to quasicontinuous domains,
\newblock \emph{Houston Journal of Mathematics}
49 (2023), no.~1, 87--103.

\bibitem{ShenEtAl}
C.~Shen, W.~K.~Ho, X.~Xi, and D.~Zhao,
\newblock Scott function spaces under one-sided FS assumptions:
counterexamples, positive results, and new directions,
\newblock arXiv preprint arXiv:2607.25832v1, 2026.

%\bibitem{XiShenZhao}
%X.~Xi, C.~Shen, and D.~Zhao,
%\newblock Sobriety of regular open algebras in second-countable $T_3$ spaces,
%\newblock \emph{Topology and its Applications}394 (2026), 109976.

\bibitem{XiXuLawson2016}
X.~Xi, L.~Xu, and J.~D.~Lawson,
\newblock Function spaces from coherent continuous domains to RB-domains,
\newblock \emph{Topology and its Applications}
206 (2016), 148--157.

\bibitem{XiYang}
X.~Xi and J.~Yang,
\newblock Coincidence of the Isbell and Scott topologies on
domain function spaces,
\newblock \emph{Topology and its Applications}
164 (2014), 197--206.

\end{thebibliography}
\end{document}